\outer\def\beginsection#1\par{\vskip0pt plus.3 \vsize\penalty-25
\vskip0pt plus-.3\vsize\bigskip\vskip\parskip
\message{#1}\leftline{\bf#1}\nobreak\medskip\noindent}
\magnification=1200
\baselineskip=15pt
\def\ni{\noindent}
\def\^{\widehat}
\def\T{{\bf T}}

\def\Z{{\bf Z}}
\def\Q{{\bf Q}}

\def\B{{\cal B}}

\def\Khat{\widehat K}
\def\hollowbox{\hbox{\vbox{\hrule\hbox{\vrule\vbox{\vbox to 8pt{\vfill\hbox to 5pt{\hfill}}}\vrule}\hrule}}}

\centerline{{\bf Traces on Noncommutative Homogeneous Spaces}}
\bigskip
\centerline{{\bf Magnus B. Landstad }}

\bigskip
The noncommutative Heisenberg manifolds constructed by Rieffel in [R2]
have turned out to provide interesting examples of $C^*$-algebras which are
similar, but not isomorphic to irrational rotation algebras as shown
by Abadie in [A2]. It was shown in [LR2] that these algebras are
special cases of a more general construction giving deformations
of $C(G/\Gamma)$ for a compact homogeneous
space $G/\Gamma$. The $C^*$-algebras obtained were denoted
$C_r^*(\^G/\Gamma,\rho)$ and their ideal structure was
determined in [LR2], in
this follow-up we shall describe the algebras more closely:
(1) The center ${\cal Z}C_r^*(\^G/\Gamma,\rho)$ is
isomorphic to  $ C(G/G_\rho^0)$ for a certain subgroup $ G_\rho^0$ of
$G$,
and (2)  there is a conditional expectation
$E^0:C_r^*(\^G/\Gamma,\rho)\mapsto C(G/G_\rho^0)$ and therefore
a 1-1 correspondence between
normalised traces  on $C_r^*(\^G/\Gamma,\rho)$ and
probability measures  on $G/G_\rho^0$. This is used to show that $C_r^*(\^G/\Gamma,\rho)$ also
 can be represented over $L^2(G/\Gamma)$ just as in the nondeformed case.

The results here generalise some of those in [A1-3] and should provide useful
tools for extending other results such as the determination
of the (ordered) K-theory and the noncommutative metrics of these algebras.

\beginsection 1. Preliminaries.

The deformations
of $C(G/\Gamma)$ constructed in [LR2] are based on the following
{\it standard
assumptions:}

 {\it (S1)  There is a compact abelian subgroup $K$ of $G$ and
a homomorphism $\rho :\Khat \mapsto G$ such that each
$\rho(s)$ commutes with $K$,

(S2) $\Gamma$ is a closed subgroup of $G$,
  each $x\in\Gamma$ commutes with $K$
and satisfies
$$\eqalign{
&B_x(s):=x\rho(s)x^{-1}\rho(-s)\in K
  \hbox{ for all } s\in  \Khat \hbox{ and }\cr
&\langle B_x(s),t\rangle =\langle B_x(t),s\rangle  \hbox{ for all }
s,t\in  \Khat,
        }$$

(S3) $G/\Gamma$ is compact and $\Gamma\cap K=\{e\}$.}

\medskip\ni
In [LR2, 4.8-11] it is then explained how one obtains a closed subgroup $K_\rho$
of $K$
by $K_\rho^\perp=\{t\in\Khat\ |\ \rho(2t)\in K\Gamma\}$
and a homomorphism $\theta:\Khat\mapsto K/K_\rho$.
The subgroup $G_\rho$  of $G$ is then defined as the closure of
$$
\{\rho(-2p)\theta(p)K_\rho\Gamma\ |\ p\in\Khat\}
$$
and it is shown  that $K_\rho\Gamma$
is a closed normal subgroup of $G_\rho$ with $G_\rho/K_\rho\Gamma$ a compact
abelian group.

\medskip
We define the partial Fourier transform by
$\widehat f(x,s):=\int_K\overline{\langle k,s\rangle }f(xk)\,dk$
for  $f\in C(G)$ as in [LR2] and
$$
C_{b,1}(G):=\Big\{f\in
C_b(G):\|f\|_{\infty,1}:=\sum_s{\sup_x|\widehat f(x,s)|}<
\infty\Big\}.
        $$
The space of functions we shall work with is
$C_1(G/\Gamma):=C(G/\Gamma)\cap C_{b,1}(G)$. With the operations
given by
$f^*(x)=\overline{f(x)}$ and
$$
 f* g(x) =
\sum_{s,t} \widehat f(x\rho(t),s)\widehat g(x\rho(-s),t)
\leqno{(1.1)}
$$
we have a Banach $^*$-algebra denoted $C_1(G/\Gamma,\rho)$.
 Its regular representation
$\mu$ over $L^2(G)$ is described in [LR2, Section 1], and
$C_r^*(\widehat G/\Gamma,\rho)$ is the $C^*$-closure of
$\mu(C_{1}(G/\Gamma,\rho))$.
Note that this definition is closely related to the Fell bundle approach in [AE1].

We refer to [LR2] for more details and other concepts not explained
here, in fact this article is unreadable without  [LR1-2].
\beginsection 2. The center of $C_r^*(\^G/\Gamma,\rho)$.

It was shown in [LR2, Theorem 4.15] that
$ C_1(G/G_\rho,\rho)$ with the product  (1.1) is a dense subalgebra of the center
${\cal Z}C_r^*(\^G/\Gamma,\rho)$.
However, this product will not be the pointwise product for functions in
$ C_1(G/G_\rho,\rho)$.
We shall see that the center actually is isomorphic to
$ C(G/G^0_\rho)$ where $G^0_\rho$ is another subgroup of $G$.
We thank the referee for discovering an error in our description,
this means that Remark 4.16 in [LR2] is incorrect.

We shall also need some other subgroups of $G$, and a guiding example in this section is to take the Heisenberg
manifolds as described in [LR2, Section 3] with $\mu$ and $\nu$
rational.
All concepts in [LR2, Section 4] are used without further explanation. For the first new  subgroup note that
$$\{\rho(-2p)\theta(p)K_\rho\Gamma\ |\ p\in\Khat\}
\supset\{\rho(-2p)\theta(p)K_\rho\Gamma\ |\ p\in K_\rho^\perp\}
=\{\widetilde\theta(-p)\theta(p)K_\rho\Gamma\ |\ p\in K_\rho^\perp\}.
        $$
The following should then be obvious:

\proclaim Lemma 2.1. Take
$K_0=\{\widetilde\theta(-t)\theta(t)K_\rho\ |\ t\in K_\rho^\perp\}$.
Then $K_\rho\Gamma\subset K_0\Gamma\subset G_\rho$ and
$$K_0^\perp=\{s\in K_\rho^\perp \ |\ \langle \theta(s),t\rangle
=\langle \theta(t),s\rangle \hbox{ for all }t\in K_\rho^\perp \}.
        $$
Furthermore,  all        $f\in C_1(G/G_\rho)$ satisfies  $\^f(x,s)=0$ for
$s\notin K_0^\perp$.

 In particular this means that we have
$\langle \theta(s),t\rangle =\langle \theta(t),s\rangle $ for all
$s,t\in K_0^\perp$,  so there is a function
$c:K_0^\perp\mapsto\T$ such that
$$\langle \theta(s),t\rangle ={c(s)c(t)\over c(s+t)} \quad
\hbox{ for all } \quad
s,t\in  K_0^\perp.
\leqno(2.1)$$

\proclaim Lemma 2.2. There is a function $b:\Khat\mapsto\T$ such that
$${b(s)c(t)\over b(s+t)}=\langle \theta(s),t\rangle
\quad\hbox{ for all } \quad
 s\in \Khat, \ t\in K_0^\perp.
\leqno(2.2)$$

 \ni{\bf Proof.}  Pick one $s_i$ from each equivalence class in
$\Khat/K_0^\perp$ and define $b$ by
$$b(s_i+t)=c(t)\langle \theta(s_i),-t\rangle    \hbox{ for  }t\in
K_0^\perp.
        $$
It is then straightforward to check that (2.2) holds.

\proclaim Lemma 2.3. For $f\in C_1(G/G_\rho)$ define
$$\Phi(f)(x)=\sum_{s\in K_0^\perp}\^f(x\rho(s),s)c(s).$$
Then $\Phi$ is a 1-1 $^*$-homomorphism  from
$ C_1(G/G_\rho,\rho)$ with the product (1.1) into $C(G)$ equipped
with the usual pointwise operations.

 \ni{\bf Proof.}
For     $f,g\in C_1(G/G_\rho)$ we have by (1.1) and the definition of
$G_\rho$ that
$$ \eqalign{
 \Phi(f* g)(x) &= \sum_{s\in K_0^\perp}(f*g) \ \^{}\
        (x\rho(s),s)c(s)\cr
&=\sum_{s,t\in K_0^\perp}\^f(x\rho(2s-t),t)\^g(x\rho(s-t),s-t)c(s)\cr
&=\sum \^f(x\rho(2s+t),t)\^g(x\rho(s),s)c(s+t)\cr
&=\sum \^f(x\rho(t)\theta(s),t)\^g(x\rho(s),s) c(s+t)\cr
&=\sum \^f(x\rho(t),t)\^g(x\rho(s),s)\langle \theta(s),t\rangle c(s+t)\cr
&=\sum \^f(x\rho(t),t)\^g(x\rho(s),s)c(s)c(t)\cr
&=\Phi(f)(x)\Phi( g)(x).\cr
        }$$
The $^*$-operation is complex conjugation in both algebras, so
$$ \eqalign{
 \Phi(f^* )(x) &= \sum_{s\in K_0^\perp}
\overline{\^f(x\rho(s),-s)}c(s)
= \sum\overline{\^f(x\rho(-s),s)}c(-s)\cr
&= \sum\overline{\^f(x\rho(s)\theta(-s),s)}c(-s)
= \sum\overline{\^f(x\rho(s),s)\langle \theta(-s),s\rangle }c(-s)\cr
&= \sum\overline{\^f(x\rho(s),s)c(s)}
=\overline{\Phi(f)(x)}. \quad\hollowbox\cr
        }$$

In order to define the subgroup $ G_\rho^0$ we need the following construction:

\proclaim Lemma 2.4. There is a continuous homomorphism $B:
G_\rho\mapsto K/ K_\rho $ satisfying $B(x)^2=e$ and $$ \langle
B(x),s\rangle= \langle x\rho(s) x^{-1} \rho(-s),s \rangle = \pm 1
\quad\hbox{ for } \quad s\in K_\rho^\perp. \leqno(2.3)$$

\ni{\bf Proof.} If $x\in K$ or $x=\rho(t)$ then
$ x\rho(s) x^{-1} \rho(-s)=e$, so we have $B(x)=e$.
It follows from (S1-S3) that for $x,y\in \Gamma$ and
$s,t\in \Khat$ that
$$\leqalignno{
B_x(s+t) &= B_x(s) B_x(t) &(2.4)\cr
B_{xy}(s) &= B_x(s) B_y(s). &(2.5)\cr
        }$$
For $x\in \Gamma$ and $ s \in K_\rho^\perp $  we proved in [LR2, Lemma 4.11] that
$B_x(2s)=e$. So if we look at the function
$s\mapsto\langle B_x(s),s\rangle $
we have for $ s,t \in K_\rho^\perp $ that
$$\eqalign{\langle B_x(s+t),s+t\rangle&= \langle B_x(s),s\rangle
\langle B_x(s),t\rangle \langle B_x(t),s\rangle \langle
B_x(t),t\rangle\cr &=\langle B_x(s),s\rangle \langle
B_x(2s),t\rangle \langle B_x(t),t\rangle\cr &=\langle
B_x(s),s\rangle
\langle B_x(t),t\rangle.
    }$$
So there is an element $B(x)\in K/ K_\rho $ such that
$$ \langle B(x),s\rangle=\langle B_x(s),s\rangle= \langle x\rho(s) x^{-1} \rho(-s),s \rangle
\quad\hbox{ for } \quad s\in K_\rho^\perp.
    $$
This holds for all
$x\in G_\rho$, and for a fixed $ s \in K_\rho^\perp $ this
expression is continuous in $x$. So (2.3-5) together with
$B_x(2s)=e$ implies that $B$ is a continuous homomorphism with
$\langle B(x),s\rangle= \pm 1$ and therefore $B(x)^2=e$ in $K/ K_\rho $.
$\quad\hollowbox$

Note the map $B$ can be defined the same way on the group $G_1$ defined
in [LR2, Lemma 4.7], but this is not needed here.

\proclaim Lemma 2.5. Define
$$G_\rho^0:=
\{yB(y^{-1})K_\rho\ |\ y\in G_\rho\}
=\{\rho(-2p)\theta(p)zB(z^{-1})K_\rho\ |\ p\in\Khat,z\in\Gamma\}^-.$$
Then the homomorphism $\Phi$ in Lemma 2.3 has dense
image in $C(G/ G_\rho^0)$.

\ni{\bf Proof.} If $f\in C_1(G/G_\rho)$, then
$ \^f(xy,s)= \^f(x,s)$  for $y\in G_\rho$.
So from the definition of $\Phi$ we have
$$\eqalign{\Phi(f)(xy)
&=\sum_{s\in K_0^\perp}\^f(xy\rho(s),s)c(s)
=\sum\^f(x\rho(s)y,s)
   \langle B_y(s),s\rangle c(s)\cr
&=\sum\^f(x\rho(s),s)
   \langle B(y),s\rangle c(s)
=\sum\^f(xB(y)\rho(s),s)c(s)\cr
&= \Phi(f)(xB(y)) .\cr
    }$$
Thus $\Phi(f)(xyB(y)^{-1}) =\Phi(f)(x) $.
$\Phi$ is a bijection between $ C_1(G/G_\rho)$
and $ C_1(G/G_\rho^0)$, so the image is dense.
$\quad\hollowbox$

We want to show that   $\Phi$ extends to an isomorphism between
${\cal Z}C_r^*(\^G/\Gamma,\rho)$ and $ C(G/G^0_\rho)$, and since $C_r^*(\^G/\Gamma,\rho)$
is defined by using the regular  representation $\mu$ this will follow from:

\proclaim Proposition 2.6. For $f\in C_1(G/G_\rho)$ the unitary
operator
$U=\sum_{s\in\Khat} b(s)L_{\rho(s)}P_s$ satisfies
$$U\mu(f)U^*=M(\Phi(f)).
        $$

\ni{\bf Proof.} From Propositions 1.3 and (1.11) in [LR2] and Lemma 2.2
above we have
$$\eqalign{U\mu(f)U^*
&=\sum_{s,t\in\Khat} b(t)L_{\rho(t)}P_t\mu(f) \overline{b(s)}L_{\rho(-s)}P_s\cr
&=\sum b(t)L_{\rho(s+t)}P_t M(f) \overline{b(s)}L_{\rho(-s-t)}P_s\cr
&=\sum b(t)L_{\rho(s+t)} M(\^f(\cdot,t-s))\overline{b(s)}L_{\rho(-s-t)}
    P_s\cr
&=\sum b(t+s)\overline{b(s)}L_{\rho(2s+t)}M(\^f(\cdot,t))L_{\rho(-2s-t)}
    P_s\cr
&=\sum_{s\in\Khat,t\in K_0^\perp}
        b(t+s)\overline{b(s)}M(\^f(\cdot\, \rho(2s+t),t))P_s\cr
&=\sum b(t+s)\overline{b(s)}\langle \theta(s),t\rangle
        M(\^f(\cdot\,\rho(t),t))P_s\cr
&=\sum c(t)M(\^f(\cdot\,\rho(t),t))P_s=M(\Phi(f)). \quad\hollowbox \cr
        }$$

Note that every $x\in G_\rho$ satisfies (S2), so by
[LR2, Theorem 4.3]
$\beta_x$ defined by $\beta_x(f)(y)=f(yx)$ extends to
a $^*$-automorphism of $C_r^*(\^G/\Gamma,\rho)$.
If we also look at part (1) and (2) of the proof of [LR2, Theorem 4.15], we see that it can be rephrased as
$${\cal Z} C_r^*(\^G/\Gamma,\rho)
=\{a\in\ C_r^*(\^G/\Gamma,\rho)\ |\ \beta_x(a)=a
\hbox{ for all } x\in G_\rho \}.\leqno{(2.6)}
$$
\proclaim Theorem 2.7. The map  $\Phi$ extends to a
$C^*$-isomorphism between
${\cal Z}C_r^*(\^G/\Gamma,\rho)$
and $ C(G/G_\rho^0)$
with the pointwise product.

\ni{\bf Proof.} This now follows, just note that in part (2) of the
proof of [LR2, Theorem 4.15] it was shown that $ \mu(C_1(G/G_\rho))$ is
dense in the center of $C_r^*(\^G/\Gamma,\rho)$.

\ni{\bf Remarks 2.8.}  The map $y\mapsto yB(y)^{-1} $ is an
isomorphism between the groups $G_\rho/K_\rho$ and $
G_\rho^0/K_\rho $. However, this will not imply  that the groups
$G_\rho$ and $ G_\rho^0 $ themselves are isomorphic or that $G/
G_\rho$ is homeomorphic to $G/ G_\rho^0$. Also note that
$G=G_\rho\iff G=G_\rho^0 $ which again by [LR2, Theorem 4.15] is
equivalent to $C_r^*(\^G/\Gamma,\rho)$ being simple. We shall see
in Section 5 that $B$ can be nontrivial and $G_\rho\neq G_\rho^0
$.
\beginsection 3. The conditional expectation onto the center and
traces on $C_r^*(\^G/\Gamma,\rho)$.

We have now proved that ${\cal Z}C_r^*(\^G/\Gamma,\rho)$ is generated by $\{\mu(f)\ |\ f\in C_1(G/G_\rho)\}$ and is via the map $\Phi$ isomorphic to
$C(G/G_\rho^0)$.
From (2.6) we get the natural conditional expectation $E$ of $C_r^*(\^G/\Gamma,\rho)$ onto its center by
$$E(a) =\int_{G_\rho/\Gamma}\beta_x(a)\,dx.
        $$
Note that $G_\rho/\Gamma$ is not
necessarily a group, but as noted in the preliminaries
$K_\rho\Gamma$
is a closed normal subgroup of $G_\rho$ with $G_\rho/K_\rho\Gamma$ a compact
abelian group.
(The same is true if we replace $G_\rho$ with $G_\rho^0 $.)
This means that the map $E$ is given by
$$E(a)=\int_{G_\rho/\Gamma}\beta_x(a)\,dx
        =\int_{G_\rho/K_\rho\Gamma}\int_{K_\rho}\beta_{yk}(a) \,dkdy.\leqno(3.1)
        $$
\proclaim Lemma 3.1. $E^0(a)=\Phi(E(a))$ defines a conditional expectation from
$C_r^*(\^G/\Gamma,\rho)$ onto $C(G/G_\rho^0)$.
For $f\in C_1(G/\Gamma,\rho)$ we have
$$E^0(\mu(f))(x)=
\sum_{s\in K_0^\perp} c(s) \int_{G_\rho^0/ K_\rho\Gamma}
\^f(xz\rho(s),s)\,dz.\leqno(3.2)
        $$

\ni{\bf Proof.}
With $\widetilde E(f)(x)=\int_{G_\rho/\Gamma}f(xz)\,dz
        =\int_{G_\rho/K_\rho\Gamma}\int_{K_\rho}f(xyk)\,dkdy $ we have
$ E(\mu(f))= \mu(\widetilde E(f))$, so
$$\eqalign{ E^0(\mu(f))(x)
    &=\Phi(\widetilde E(f))(x)\cr
    &=\sum_{s\in K_0^\perp} c(s) \widetilde E(f) \ \^{}\ (x\rho(s),s)\cr
    &=\sum c(s)\int_{G_\rho/\Gamma}\^f(x\rho(s)z,s)\,dz.\cr}$$
Since $ s\in K_0^\perp \subset K_\rho^\perp $ we have
$$\eqalign{ E^0(\mu(f))(x)
&=\sum c(s)\int_{G_\rho/K_\rho\Gamma}\^f(x\rho(s)y,s)\,dy\cr
&=\sum c(s)\int \^f(xy\rho(s)B_{y^{-1}}(s),s)\,dy\cr
&=\sum c(s)\int \^f(xy\rho(s),s) \langle B_{y^{-1}}(s) ,s\rangle \,dy.\cr}
    $$
From Lemma 2.4 we have
$\langle B_{y^{-1}}(s) ,s\rangle =\langle B(y^{-1}) ,s\rangle $, so
$$\eqalign{ E^0(\mu(f))(x)
&=\sum c(s)\int_{G_\rho/K_\rho\Gamma} \^f(xy B(y^{-1}) \rho(s),s) \,dy\cr
&=\sum c(s)\int_{G_\rho^0/K_\rho\Gamma} \^f(xz\rho(s),s) \,dz.\cr}
    $$
\proclaim Lemma 3.2. For $f,g\in C_1(G/\Gamma,\rho)$ we have
$E^0(\mu(f*g))= E^0(\mu(g*f))$.

\ni{\bf Proof.} There are no problems in interchanging integrals and sums here, so for
$s\in K_0^\perp$
$$\eqalign{
\int_{G_\rho^0/K_\rho\Gamma} (f*g) \ \^{}\ (xz\rho(s),s) \,dz
&=\int\sum_{t\in\Khat}\^f(xz\rho(2s-t),t)\^g(xz\rho(s-t),s-t)\,dz\cr
&= \sum\int \^f(xz\rho(s+t),s-t)\^g(xz\rho(t),t)\,dz.\cr
    }$$
Here we used the substitution $t\mapsto s-t$,
we continue with the substitution
$z K_\rho \mapsto z\rho(-2t)\theta(t)K_\rho$
to get
$$\int (f*g) \ \^{}\ (xz\rho(s),s) \,dz
=\sum\int \^f(xz\rho(s-t),s-t)\^g(xz\rho(-t),t)
\langle \theta(t) ,s-t+t \rangle \,dz.
    $$
From (4.9) in [LR2] we have that
$\widetilde\theta(s)\Gamma=\rho(2s)\Gamma$ for $s\in K_\rho^\perp$,
which together with
$\langle \theta(t) ,s \rangle= \langle \widetilde\theta(s) ,t \rangle $ gives
$$\eqalign{\int (f*g) \ \^{}\ (xz\rho(s),s) \,dz
    &=\sum\int \^f(xz\rho(s-t),s-t) \^g(xz\rho(2s-t),t) \,dz\cr
    &=\int (g*f) \ \^{}\ (xz\rho(s),s) \,dz.\cr
    }$$
This holds for all $s\in K_0^\perp$, so from (3.2) we have
$E^0(\mu(f*g))= E^0(\mu(g*f))$. $\quad\hollowbox$

We can now prove the following generalisation of [A3, Corollary 3.11]:
\proclaim Theorem 3.3. The conditional expectation from
$C_r^*(\^G/\Gamma,\rho)$ onto $C(G/G_\rho^0)$ given by
$E^0=\Phi\circ E$ satisfies $\alpha_x\circ E^0=E^0\circ\alpha_x $
and $E^0(ab)=E^0(ba)$. There is a 1-1 correspondence between
normalised traces $\tau$ on $C_r^*(\^G/\Gamma,\rho)$ and
probability measures $\nu$ on $G/G_\rho^0$ given by
$\tau=\nu\circ E^0$.
$\tau$ is faithful if and only if
$\nu$ is.

\ni{\bf Proof.} The $\alpha$-invariance is obvious.
Since
$\mu(C_1(G/\Gamma,\rho))$ is dense in $C_r^*(\^G/\Gamma,\rho)$, the
first part follows from Lemma 3.2. Hence
$\tau=\nu\circ E^0$ is a normalised trace for all probability measures
$\nu$ on $G/G_\rho^0$. Since $E^0$ is faithful, it
is immediate that $\tau$ is faithful if and only if $\nu$ is.

Conversely, if  $\tau$ is a normalised trace on
$C_r^*(\^G/\Gamma,\rho)$, it follows from [LR2, Lemma 4.12] that
$\tau\circ\beta_k=\tau$ for $k\in K_\rho$.
So $\tau(a)=\tau(b)$ where
$b=\int_{K_\rho}\beta_k(a)\,dk$.
From [LR2, Lemma 4.14] it follows that
$\tau\circ\beta_y(b)=\tau(b)$ for $y\in G_\rho$, so
$$\tau(a)
=\int_{G_\rho/K_\rho\Gamma}\int_{K_\rho}\tau\circ\beta_{yk}(a)\,dkdy
=\tau(E(a)).
$$
Let $\nu$ be the measure on $G/G_\rho^0$
given by $\nu(f)=\tau(\Phi^{-1}(f))$ for $f\in C(G/G_\rho^0)$,
then $\tau(a)=\tau (E(a)) =\nu(\Phi\circ E(a)) =\nu(E^0(a))$.
$\quad\hollowbox$

From [LR2, Theorem 4.15] we now have
\proclaim Corollary 3.4. If $G=G_\rho$ (which is equivalent to $G=G_\rho^0$),
there is a unique trace on the
simple $C^*$-algebra
$C_r^*(\^G/\Gamma,\rho)$.

\beginsection 4. Quasi-invariant measures on $G/\Gamma$ and
representations over $L^2(G/\Gamma)$.

In this section we shall look at traces on $C_r^*(\^G/\Gamma,\rho)$
obtained from a $G$-quasi-invariant measure $\nu$ on $G/G_\rho^0$.
We shall see that the  corresponding GNS-representation can be realised
over $L^2(G/\Gamma,\nu)$.

If $\nu$ is a $G$-quasi-invariant measure on $G/G_\rho^0$, there is a
function $h\in C(G\times G/G_\rho^0)$ such that
$\nu(\alpha_x(f))=\nu(h(x,\cdot)f)$ for $f\in C(G/G_\rho^0)$, in fact
$h(x,y)={r(xy)\over r(y)}$ where $r$ is a continuous rho-function on
$G$ corresponding to $\nu$.
If we  take
$z_x=\Phi^{-1}[h(x,\cdot)]$
and  use Theorem 3.3 on such measures, we get
\proclaim Corollary 4.1. If $\nu$ is a $G$-quasi-invariant
probability measure on $G/G_\rho^0$ and $\tau=\nu\circ E^0$, there is
a continuous function
$x\in G\mapsto z_x\in {\cal Z}C_r^*(\^G/\Gamma,\rho)^+$
 such that
$\tau\circ\alpha_x(a)=\tau(z_xa)$ for all
$a\in C_r^*(\^G/\Gamma,\rho)$.

Since $G_\rho^0/\Gamma$ has a $G_\rho^0$-invariant probability measure
and $G/\Gamma$ is compact, it follows
that $G/\Gamma$ has a $G$-invariant probability
measure if and only if $G/G_\rho^0$ has one. And if so, it is
unique. However, note that the
compactness of $G/\Gamma$ does {\it not} imply the existence of
a $G$-invariant probability measure.
\proclaim Corollary 4.2. If $G/\Gamma$ has a $G$-invariant probability
measure, then there is a unique normalised
$\alpha_G$-invariant trace on $C_r^*(\^G/\Gamma,\rho)$.

The regular representation of $C_r^*(\^G/\Gamma,\rho)$ is over
$L^2(G)$, but it seems natural to also represent it over
$L^2(G/\Gamma)$. This is obtained by using the GNS-representation
obtained from
 $\tau$ as in  Corollary 4.1.
\proclaim Lemma 4.3. $G/\Gamma$ has a
$G$-quasi-invariant probability measure $\nu$  such that
$\nu(g^**f)=\nu(\overline gf) $ for $f,g\in C_1(G/\Gamma)$.

\ni{\bf Proof.} First we shall need the closed subgroup
$G_2=\{\rho(t)K\Gamma\ |\ t\in\Khat\}^-$ of $G$. Exactly as in
[LR2, Lemma 4.7] one proves that $K\Gamma$
is a closed, normal subgroup of $G_2$, and $G_2/K\Gamma$ is a compact
abelian group. If we now take a
$G$-quasi-invariant probability measure  on $G/G_2$, Haar-measures on
$G_2/K\Gamma$ and $K$, then
$$\int f\ d\nu
=\int_{G/G_2}\int_{G_2/K\Gamma}\int_{K}f(xyk)\,dk\,dy\,dx
        $$
defines a $G$-quasi-invariant probability measure on $G/\Gamma$.
So for  $f,g\in C_1(G/\Gamma)$ we have
$$\eqalign{\int_{K}g^**f(xk)\,dk
&=\int_{K}\sum_{s,t}\overline{\^g(xk\rho(t),-s)}\^f(xk\rho(-s),t)\,dk\cr
&=\sum_{t}\overline{\^g(x\rho(t),t)}\^f(x\rho(t),t).\cr
\nu(g^**f)
&=\int_{G/G_2}\int_{G_2/K\Gamma}\sum_{t}\overline{\^g(xy\rho(t),t)}
  \^f(xy\rho(t),t)\,dy\,dx\quad (y\mapsto y\rho(-t))\cr
&=\int_{G/G_2}\int_{G_2/K\Gamma}\sum_{t}\overline{\^g(xy,t)}
  \^f(xy,t)\,dy\,dx\cr
&=\int_{G/G_2}\int_{G_2/K\Gamma}\int_K\overline{g(xyk)}
  f(xyk)\,dk\,dy\,dx\cr
&=\nu(\overline gf).
        }$$
\proclaim Lemma 4.4. Let  $\nu$ be a
$G$-quasi-invariant probability measure on $G/\Gamma$. Then there is a
function $\phi\in C_c(G)$ such that for all
$f\in C_1(G/\Gamma)$
$$\langle\mu(f)\phi\, |\, \phi\rangle=\int_{G/\Gamma} f\,d\nu.\leqno(4.1)
$$

\ni{\bf Proof.}
Let $r$ be a continuous rho-function as in
[FD, Section III.13.2] or [LR1, Section 2], we may assume that
$r(xk)=r(x)$ for $k\in K$. So for $f\in C_c(G)$
$$
\int_{G/\Gamma}\int_\Gamma f(xh)\,dh\,d\nu(x)=\int_G
f(x)r(x)\,dx .
        $$
Since $G/\Gamma$ is compact, there is a function $\phi_0\in C_c(G)$
such that
$$\int_\Gamma\phi_0(xh)\,dh=1\quad\hbox{for all}\quad x\in G.\leqno(4.2)
        $$
$K$ is compact and commutes with $\Gamma$, so we may assume that
$\phi_0(xk)=\phi_0(x)$ for $k\in K$.
Now take
$\phi(x)=[\phi_0(x^{-1})r(x^{-1})\Delta(x^{-1})]^{1\over2}$.
Since $P_0\phi=\phi$, we get from [LR2, Proposition 1.3] that
$$\eqalign{\langle\mu(f)\phi\, |\, \phi\rangle
  &=\langle M(f)\phi\, |\, \phi\rangle=\int_Gf(x^{-1})\phi(x)^2\,dx\cr
  &=\int_Gf(x^{-1})\phi_0(x^{-1})r(x^{-1})\Delta(x^{-1})\,dx
  =\int_Gf(x)\phi_0(x)r(x)\,dx\cr
  &=\int_{G/\Gamma}\int_\Gamma f(y)
    \phi_0(yh)\,dh\,d\nu(y)
   =\int_{G/\Gamma} f(y)\,d\nu(y).\quad \hollowbox\cr
        }$$

These two lemmas show that  $C_r^*(\^G/\Gamma,\rho)$ -- which was
defined over $L^2(G)$ -- also can be represented over
$L^2(G/\Gamma,\nu)$. This is done as follows: Define
$V:C_1(G/\Gamma)\mapsto L^2(G)$  by $Vf=\mu(f)\phi$. Then $\langle
Vf\, |\,Vg\rangle=\langle\mu(g^**f)\phi\, |\, \phi\rangle
=\nu(g^**f)=\nu(\overline gf)$, so $V$ extends to a partial
isometry from $L^2(G/\Gamma,\nu)$ into $L^2(G)$ with  $V^*V=I$. We
also have
 $$V^*\mu(f)Vg=V^*\mu(f)\mu(g)\phi=V^*\mu(f*g)\phi=V^*Vf*g=f*g.
        $$
Thus the GNS-representation of
$C_r^*(\^G/\Gamma,\rho)$ corresponding to $\tau$ is really over
$L^2(G/\Gamma,\nu)$ and
$V$ sets up an equivalence with a sub-representation of the
regular representation $\mu$. It is faithful, since
$\langle V^*aV 1\, |\, 1\rangle =\langle a\phi\, |\, \phi\rangle =\tau(a)$
and $\tau$ is faithful. Thus we have shown:
\proclaim Theorem 4.5.  If $\nu$ is a
$G$-quasi-invariant probability measure on $G/\Gamma$ as in Lemma 4.3,
then the GNS-representation of
$C_r^*(\^G/\Gamma,\rho)$ corresponding to  $\nu$ is a faithful representation
over $L^2(G/\Gamma,\nu)$ and
equivalent to a sub-representation of the
regular representation $\mu$.

\beginsection 5. The Heisenberg manifolds revisited.

Let us go back to the Heisenberg manifolds in [R2],
we shall use our description in [LR2, Section 3];
so $C_r^*(\^G/\Gamma,\rho)\cong D_{\mu,\nu}$ in the
terminology of [LR2] and [A1-3], we only look at the case with $c=1$.
We only state the results and leave the computations to the reader.

If $\mu,\nu\in\Q$ one finds that $ K_\rho^\perp = q\Z$, where
$q$ is the smallest integer $\neq 0$ with both $2\mu q$ and $2\nu
q\in\Z$. $\widetilde\theta(t)=\theta(t)=(-1)^{4q\mu\nu t}$ for $t
\in q\Z$, so $ K_\rho^\perp = K_0^\perp $. For any $g=(x,y,z)\in
G$ define $B(g)=e^{2\pi iq(\mu y-\nu x)}$ and $B(g)$ satisfies
(2.3) for $g\in\Gamma$. If e.g.\  $\mu=\nu={1\over2}$, then $q=1$
and $ K_\rho =\{e\}$, but $B$ is not trivial on $\Gamma$, so $
G_\rho\neq  G_\rho^0 $. However, since $B$ is defined on the whole
group $G$, we have $G_\rho\cong G_\rho^0 $ and $G/ G_\rho$ is
homeomorphic to $G/ G_\rho^0$.

On the other hand if
$\mu$ or $\nu$ is irrational, then
$ K_\rho = K_0=K$, so
$ G_\rho = G_\rho^0 =\{\rho(2n)K\Gamma\ |\ n\in\Z\}\bar{\ }\bar{\ }$ 
and is a normal subgroup of $G$ with
$G/G_\rho\cong\T^2/T_0$ where $T_0$ is the closed
subgroup generated by $(\exp{4\pi i\mu},\exp{4\pi i\nu})$.
So Theorem 3.3 is our version of [A3, Corollary 3.11].

Many of the structural results about $ D_{\mu,\nu}$ in [A1-3]
can be proved using the present description.
Let us briefly illustrate this by showing the
existence of projections \'a la Rieffel
using only functions in $ C_1(G/\Gamma)$.
Let
$$\eqalign{
 f(x,y,z)&= F_\mu (x)\qquad
g(x,y,z)= G_\mu (x)z \exp(-2\pi i[x]y)\cr
h(x,y,z)&= F_\nu (y)\qquad
k(x,y,z)= G_\nu (y)z \exp(-2\pi ix([y]-y))
    }$$
where $F_\mu$ and $G_\mu$ are continuous functions satisfying a
slightly modified version of [R1, Theorem 1.1 (1-3)]. In
particular we need $G_\mu(0)=0$ in order to make $g$ and $k$
continuous. Computations as in [LR2, Proposition 3.6] show that
$p=f+g+g^*$ and  $q=h+k+k^*$ are projections in $D_{\mu,\nu}$ with
$\tau(p)=2\mu$ and $\tau(q)=2\nu$ for any normalised trace $\tau$.

This is used in [A1, Theorem 1] to show that $\Z+2\mu \Z+2\nu\Z\subset K_0(D_{\mu,\nu}) $.
To show equality Abadie proves in [A3] that  $D_{\mu,\nu}$ can be embedded in an AF-algebra.
This can also be done using only functions in $ C_1(G/\Gamma)$. The function
$$ w(x,y,z)= z \exp(-2\pi i[x]y)
    $$
is in $ L^\infty(G/\Gamma)$, but is not continuous. Since
$w=\^w(\cdot,1)$, the regular representation in [LR2, Proposition
1.3] can still be used on $w$ to get a unitary operator $W$. The
$C^*$-algebra $\B$ generated by $D_{\mu,\nu}$ and $W$ is invariant
under the automorphisms $\beta_t $ and we have $\beta_t(W)=tW$ for
all $t\in\T$. It is then standard to show that  $\B$ is in fact
the crossed product of a norm-closed subalgebra of $
L^\infty(\T^2)$ with $\Z$ as in [A3, Theorem 2.3], and by
classical results by Pimsner in [P] it follows that  $\B$ can be
embedded in an AF-algebra. Abadie uses this to determine the
ordered K-theory of $D_{\mu,\nu}$ and to describe when two such
algebras are isomorphic: In most cases $D_{\mu,\nu}\cong
D_{\mu',\nu'}$ if and only if $(\mu,\nu)$ and $(\mu',\nu')$ belong
to the same orbit under the natural action of GL(2,$\Z$) on
$\T^2$, see [AE2, Theorem 2.2] and [A3, Corollary 3.17].

Note that the functions above (except $w$) can be taken to be
$C^\infty $-functions, so also the cyclic cohomology of
$D_{\mu,\nu}$ can  be studied. It should then be possible to
extend the results for the Heisenberg manifolds to the more
general algebras $C_r^*(\^G/\Gamma,\rho)$ described here by
finding functions in $C_{1}(G/\Gamma,\rho)$ having the right
properties. It is our belief that the presentation of these
algebras given in [LR1-2] and here will be useful for such
constructions. The noncommutative metrics studied by Rieffel and
Weaver in [R3] and [W] are examples of this.

\beginsection References.

\baselineskip=12pt
\frenchspacing
\item{[A1]} B.~Abadie, "Vector bundles" over quantum Heisenberg manifolds.
{\it Algebraic methods in operator theory,}
307--315, Birkhäuser Boston, Boston, (1994).
\smallskip
\item{[A2]} B.~Abadie, Generalized fixed-point algebras of certain
actions on crossed products,
{\sl  Pacific J.~Math. {\bf 171}}  (1995), 1--21.
\smallskip
\item{[A3]} B.~Abadie,  The range of traces on quantum Heisenberg manifolds,
{\sl Trans. Amer. Math. Soc. {\bf 352}} (2000), 5767--5780 (electronic).
\item{[AE1]} B.~Abadie and R.~Exel, Deformation quantization via Fell Bundles, {\sl Math. Scand.} (to appear).
\smallskip
\item{[AE2]} B.~Abadie and R.~Exel, Hilbert $C\sp *$-bimodules over commutative $C\sp *$-algebras and an
isomorphism condition for quantum Heisenberg manifolds.
{\sl Rev. Math. Phys. {\bf 9}} (1997), 411--423.
\smallskip
\item{[FD]}   J.~M.~G.~Fell and R.~Doran,
  {\it Representations of $^*$-algebras, locally compact groups, and
Banach $^*$-algebraic bundles}, Academic Press, (1988).
\smallskip
\item{[LR1]} M.~B.~Landstad and I.~Raeburn, Twisted dual-group
algebras: Equivariant deformations of $C_0(G)$,
{\sl J.~Funct. Anal. {\bf 132}}  (1995), 43--85.
\smallskip
\item{[LR2]} M.~B.~Landstad and I. Raeburn,
 Equivariant deformations of homogeneous spaces,
{\sl J.~Funct. Anal. {\bf 148}}  (1997), 480--507.
\smallskip
\item{[P]} M.~V.~Pimsner, Embedding some transformation group $C^*$-algebras into AF-algebras,
{\sl Ergodic Theory Dynamical Systems {\bf 3}}  (1983), 613--626.
\smallskip
\item{[R1]} M.~A.~Rieffel, $C^*$-algebras associated with irrational rotations.
{\sl  Pacific J.~Math. {\bf 91}}  (1981), 415--429.
\smallskip
\item{[R2]} M.~A.~Rieffel, Deformation quantization of Heisenberg
manifolds, {\sl Comm. Math. Phys. {\bf 122}} (1989), 531--562.
\smallskip
\item{[R3]} M.~A.~Rieffel, Metrics on states from actions of compact groups.
{\sl Doc. Math. {\bf 3}} (1998),
215--229 (electronic).
\smallskip
\item{[W]} N.~Weaver, Sub-Riemannian metrics for quantum Heisenberg manifolds,
{\sl J. Operator Theory {\bf 43}}
(2000), 223--242.
\bigskip
Department of Mathematical Sciences

The Norwegian University of Science and Technology

N-7491 Trondheim

Norway

\end